\documentclass[12pt,leqno]{article}
\usepackage{amsmath,amsthm}

\def\init{\setcounter{equation}{0}}
\newtheorem{theorem}{Theorem}[section]

\newcommand{\R}{{\bf R}}

\newtheorem{lemma}{Lemma}[section]

\newcommand{\e}{{\varepsilon}}
\newcommand{\D}{\mathcal{D}}

\title{A new approach to  hyperbolic inverse problems.
\author{G.Eskin, \ \ \  Department of Mathematics, UCLA,\\ Los Angeles,
CA 90095-1555, USA. \ E-mail: eskin@math.ucla.edu}
}

\begin{document}

\maketitle
\begin{abstract}
We present a modification of the BC-method 
in the inverse hyperbolic problems.
The main novelty is the study of the restrictions of the solutions to
the characteristic surfaces instead of the fixed time hyperplanes. 
The main result is that the time-dependent Dirichlet-to-Neumann operator
prescribed on a part of the boundary uniquely determines the coefficients
of the self-adjoint hyperbolic operator up to a diffeomorphism and
a gauge transformation.  In this paper we prove the crucial local step.
The global step of the proof will be presented in the forthcoming paper.
\end{abstract}

\section{Introduction.}
\label{section 1}
\init

Consider a hyperbolic equation of the form 
\begin{eqnarray}                               \label{eq:1.1}
Lu\stackrel{def}{=}
\frac{\partial^2 u}{\partial t^2}
+\sum_{j,k=1}^n\frac{1}{\sqrt{g(x)}}
\left(-i\frac{\partial}{\partial x_j}+A_j(x)
\right)
 \sqrt{g(x)}g^{jk}(x)
\left(-i\frac{\partial}{\partial x_k}+A_k(x)\right)u
\nonumber
\\
+V(x)u=0
\end{eqnarray}
in $\Omega\times(0,T_0)$,  where  $\Omega$ is a smooth bounded domain
in $\R^n$, all coefficients in (\ref{eq:1.1}) are 
$C^\infty(\overline{\Omega})$ functions, 
$ \|g^{jk}(x)\|^{-1}$ is the 
metric tensor
 in $\overline{\Omega}$,
$g(x)=\det\|g^{jk}\|^{-1}$.  We assume that
\begin{equation}                                 \label{eq:1.2}
u(x,0)=u_t(x,0)=0 \ \ \mbox{in}\ \ \ \Omega,\ \ \
u\left|_{\partial\Omega\times(0,T_0)}\right. = f(x,t).
\end{equation}

We shall study the inverse boundary problem with boundary data given
on the part of the boundary.  Let $\Gamma_0$ 
be an open subset of $\partial\Omega$ ( in particular, $\Gamma_0=
\partial\Omega$ )
and let the D-to-N (Dirichlet-to-Neumann) operator
\begin{equation}                                   \label{eq:1.3}
\Lambda f=\sum_{j,k=1}^n g^{jk}(x)\left(\frac{\partial u}{\partial x_j}
+iA_j(x)u\right)\nu_k\left
(\sum_{p,r=1}^ng^{pr}(x)\nu_p\nu_r\right)^{-\frac{1}{2}}
\left|_{\Gamma_0\times(0,T_0)}\right.
\end{equation}
be given for all $f$ with the supports in $\Gamma_0\times (0,T_0]$.
Here $\nu=(\nu_1,...,\nu_n)$ is the unit exterior normal to 
$\partial\Omega$ with respect to the Euclidian metric.
If $F(x)=0$ is the equation of $\partial\Omega$ in some neighborhood
then (\ref{eq:1.3}) has the following form in this neighborhood:
 \begin{equation}                                   \label{eq:1.4}
\Lambda f=\sum_{j,k=1}^n g^{jk}(x)\left(\frac{\partial u}{\partial x_j}
+iA_j(x)u\right)F_{x_j}(x)
\left(\sum_{p,r=1}^ng^{pr}(x)F_{x_p}F_{x_r}
\right)^{-\frac{1}{2}}
\left|_{F(x)=0,0<t<T_0}\right. .
\end{equation}
Let
\begin{equation}                                   \label{eq:1.5}
y=y(x)
\end{equation}
be a $C^\infty$ diffeomorphism  of $\overline{\Omega}$ 
onto $\overline{\Omega_0}$ 
such that
$\Gamma_0\subset \partial\Omega_0$, the Jacobian $\det\frac{\D y}{\D x}
\neq 0$ in $\overline{\Omega}$ and 
\begin{equation}                                \label{eq:1.6}
y=x \ \ \ \mbox{on}\ \ \ \Gamma_0.
\end{equation}
If we change variables $y=y(x)$ in (\ref{eq:1.1}) 
we get an equation  in $\Omega_0$ of the same form as (\ref{eq:1.1}):
\begin{eqnarray}                               \label{eq:1.7}
\frac{\partial^2 v}{\partial t^2}
+\sum_{j,k=1}^n
\frac{1}{\sqrt{g_0(y)}}
\left(-i\frac{\partial}{\partial y_j}
+A_j^{(0)}(y)\right)
\sqrt{g_0(y)}
g_0^{jk}(y)
\left(-i\frac{\partial}{\partial y_k}
+   A_k^{(0)}(y)\right)v(y,t)
\nonumber
\\
+V^{(0)}(y)v(y,t)=0,
\ \ \ \ \ \ \ \ \ \ \ \ \ \ \ \ \ \ 
\end{eqnarray}
where 
$v(y(x),t)=u(x,t)$,
\begin{equation}                                 \label{eq:1.8}
\|g_0^{jk}(y(x))\|=\left(\frac{\D y}{\D x}\right)\|g^{jk}(x)\|
\left(\frac{\D y}{\D x}\right)^T,
\end{equation}
\begin{equation}                                  \label{eq:1.9}
g_0(y)=\det \|g_0^{jk}(y)\|^{-1},
\ \ \ \ V^{(0)}(y(x))=V(x)
\end{equation}
and
\begin{equation}                                   \label{eq:1.10}
y^*(A^{(0)})=A,
\end{equation}
where
\begin{equation}                                   \label{eq:1.11}
A=\sum_{j=1}^nA_j(x)dx_j,\ \ \ 
A^{(0)}=\sum_{j=1}^nA_j^{(0)}(y)dy_j.
\end{equation}
The D-to-N operator
$\Lambda^{(0)}$ corresponding to (\ref{eq:1.7}) has the form 
\begin{equation}                                    \label{eq:1.12}
\Lambda^{(0)}f=\sum_{j,k=1}^n g_0^{jk}(y)
\left(\frac{\partial v}{\partial y_j} + iA_j^{(0)}(y)v\right)\nu_k
\left(\sum_{p,r=1}^ng_0^{pr}(y)\nu_p\nu_r
\right)^{-\frac{1}{2}}
\left|_{\Gamma_0\times (0,T_0)}\right.,
\end{equation}
where $v(y,t)\left|_{\Gamma_0\times(0,T_0)}\right.
=f(y,t)$  since (\ref{eq:1.5}) is the identity on $\Gamma_0$.

It follows from (\ref{eq:1.4}) and (\ref{eq:1.6}) that
\begin{equation}                                   \label{eq:1.13}
\Lambda^{(0)}f=\Lambda f, \ \ \forall f,\ \  \mbox{supp\ }f \subset
\Gamma_0\times (0,T_0].
\end{equation}

Denote by $G_0(\overline{\Omega})$ the group of $C^\infty(\overline{\Omega})$
functions such that $c(x)\neq 0$ on $\overline{\Omega}$ and $c(x)=1$ on
$\Gamma_0$.  We say that $A(x)=(A_1(x),...,A_n(x))$ and 
$A^{(1)}=(A_1^{(1)}(x),...,A_n^{(1)}(x))$ are gauge equivalent in $\Omega$
if there exists $c(x)\in G_0(\overline{\Omega})$ such that
\[
A_j^{(1)}(x)=A_j(x)-ic^{-1}(x)\frac{\partial c}{\partial x_j},\ \ 
1\leq j\leq n.
\]
If $u^{(1)}(x,t)=c^{-1}(x)u(x,t)$ where $u(x,t)$ is the solution of 
(\ref{eq:1.1}) then $L^{(1)}u^{(1)}=0$ where $L^{(1)}$ is the same
as $L$ with $A(x)$ replaced by $A^{(1)}(x)$.  We shall
write for the brevity that $c\circ L^{(1)}=L$.

Denote \[
T_*=\sup_{x\in\overline{\Omega}}d(x,\Gamma_0),
\]
 where 
$d(x,\Gamma_0)$ is the distance in $\overline{\Omega}$ with respect to metric
tensor
$\|g^{jk}\|^{-1}$ from $x\in\Omega$ to $\Gamma_0$.
\begin{theorem}                                \label{theo:1.1}
Suppose $T_0>2T_*$ and $L$ and $L_0$ are
 operators 
of the form (\ref{eq:1.1}) and (\ref{eq:1.7}) in
$\Omega$  and $\Omega_0$ respectively,
where $A_j(x),V(x)$ and $A_j^{(0)},V^{(0)}$ are real valued, $1\leq j\leq n$.  
Let $\Lambda$ be the D-to-N operator corresponding
to (\ref{eq:1.1}).  If $\Lambda^{(0)}$ is the D-to-N corresponding
to (\ref{eq:1.7}) and $\Lambda=\Lambda^{(0)}$ on $\Gamma_0\times(0,T_0)$ then
there exists a diffeomorphism (\ref{eq:1.5}) 
and a gauge transformation $c(x)\in G_0(\overline{\Omega})$ such that
\[
c\circ y^{-1}\circ L^{(0)}=L.
\]
\end{theorem}

The first result in this direction was obtained in [I].  
The most general results were obtained by the BC-method 
(see [B1], [B2], [K], [KK] and {[KKL]).  

An important class of inverse problems with the data
given on a part of the boundary is the inverse problems in
domains with obstacles.
 In this case $\Omega$ is a ball $B_R$ with removed subdomains 
$\Omega_1,...,\Omega_m$,  called obstacles.
The D-to-N operator $\Lambda$ is given
on $\partial B_R\times(0,T_0)$ and the zero Dirichlet boundary
conditions are satisfied on $\partial\Omega_j\times(0,T_0),\ j=1,...,m$.
In  [E1] such inverse problems were considered
in the connection with the Aharonov-Bohm effect assuming that
$B_R$ contains one or several convex obstacles.

In the present paper we developed a modification of the BC-method.
The main novelty is the study of the restrictions of the solutions
to the characteristic surface instead of the restrictions to
the hyperplane $t=const$ as in BC-method.
The proof of Theorem \ref{theo:1.1}
consists of two steps.  In the first step we 
prove that knowing $\Lambda$ in the neighborhood of 
$\Gamma\times(0,T)$, where 
$\Gamma\subset \Gamma_0,\ 0<T<T_0,$ one can recover the coefficients 
of the equation (\ref{eq:1.1}) in some neighborhood of $\Gamma$
up to a diffeomorphism and a gauge transformation.
This
will be done in \S2.  In \S3 we prove some lemmas used in \S 2.
The global step of the proof of Theorem \ref{theo:1.1}  will be given
in the forthcoming paper.
A generalization to the case of Yang-Mills potentials is considered in [E4].

\section{The local step.}
\label{section 2}
\init

Let $\Gamma$ be an open subset of $\Gamma_0$ 
 and let $U_0$ be a  neighborhood of $\Gamma$.
Let $(x',x_n)$ be a system of coordinates in $U_0$ 
 such that the equation of $\Gamma_0$ is 
$x_n=0$ and $x'=(x_1,...,x_{n-1})$ are coordinates on 
$\Gamma_0\cap U_0$.

We introduce semi-geodesic coordinates for $L$.
Denote by $\varphi_{n}(x)$ the solution of the equation
(c.f.[E3])
\begin{eqnarray}                                    \label{eq:2.1}
\sum_{j,k=1}^n g^{jk}(x)\frac{\partial\varphi_{n}}{\partial x_j}
\frac{\partial\varphi_{n}}{\partial x_k}=1,\ \ 0\leq x_n<\delta,
\\
\varphi_{n}(x',0)=0.
\nonumber
\end{eqnarray}
Also denote by $\varphi_{p}(x),\ 1\leq p\leq n-1$, the solutions
of the equations 
\begin{eqnarray}                                    \label{eq:2.2}
\sum_{j,k=1}^n g^{jk}(x)\frac{\partial\varphi_{n}}{\partial x_j}
\frac{\partial\varphi_{p}}{\partial x_k}=0,
\\
\varphi_{p}(x',0)=x_p, \ \ 1\leq p\leq n-1.
\nonumber
\end{eqnarray}
We have that $y=\varphi(x)=(\varphi_{1}(x),...,\varphi_{n}(x))$ 
exists for
$0\leq x_n\leq \delta$ when $\delta$ is small and the Jacobian
\[
J(x)=\left|\frac{\mathcal{D}\varphi}{\mathcal{D}x}\right|\neq 0
\ \ \mbox{in} \ \ \overline{\Gamma}\times [0,\delta].
\]
Note that $\varphi(x',0)=x'.$
Changing variables $y=\varphi(x)$  in (\ref{eq:1.1}) we get 
\[
\hat{L}\hat{u}=0,
\]
where $\hat{u}(y,t)=u(x,t)$ and
\begin{eqnarray}                             \label{eq:2.3}
 \ \ \ \ 
\\
\hat{L} = \frac{\partial^2}{\partial t^2}
+\sum_{j,k=1}^n\frac{1}{\sqrt{\hat{g}(y)}}\left(-i\frac{\partial}{\partial y_j}
+\hat{A}_j(y)\right)
\sqrt{\hat{g}(y)}\hat{g}^{jk}(y)
\left(-i\frac{\partial}{\partial y_k}+\hat{A}_k(y)\right)
\nonumber
\\
+\hat{V}(y),
\nonumber
\end{eqnarray}
\begin{equation}                               \label{eq:2.4}
\hat{g}^{jk}(\varphi(x))=\sum_{p,r=1}^n g^{pr}(x)
\frac{\partial\varphi_{j}}{\partial x^p}
\frac{\partial\varphi_{k}}{\partial x^r},
\ \ \hat{g}(y)=\det\|\hat{g}^{jk}(y)\|^{-1},
\ \ \hat{V}(\varphi(x))=V(x) 
\end{equation}
and 
\[
A_k(x)=
\sum_{j=1}^n\hat{A}_j(\varphi(x))\frac{\partial\varphi_j(x)}{\partial x_k}.
\]

It follows from (\ref{eq:2.1}), (\ref{eq:2.2}) that
\begin{equation}                             \label{eq:2.5}
\hat{g}^{nn}(y)=1,\ \ \hat{g}^{nj}=\hat{g}^{jn}=0,\ \
1\leq j\leq n-1.
\end{equation}

Let
\begin{equation}                              \label{eq:2.6}
A_j'(y)=-\frac{i}{2}(\sqrt{\hat{g}})^{-1}\frac{\partial\sqrt{\hat{g}}}
{\partial y_j}
=-\frac{i\hat{g}_{y_j}}{4\hat{g}},\ \ 1\leq j\leq n.
\end{equation}
Then we can rewrite (\ref{eq:2.3}) in the form
\begin{eqnarray}                                \label{eq:2.7}
\hat{L}\hat{u}=\frac{\partial^2\hat{u}}{\partial t^2}
+\left(-i\frac{\partial}{\partial y_n}+\hat{A}_n(y)
+A_n'(y) \right)^2\hat{u}
\ \ \ \ \ \ \ \ \ \ \ \ \ \ \ 
\\
+\sum_{j,k=1}^{n-1}\left(-i\frac{\partial}{\partial y_j}+\hat{A}_j
+A_j'(y)\right)
\hat{g}^{jk}\left(-i\frac{\partial}{\partial y_k} +
\hat{A}_k+A_k'\right)\hat{u}
+V_1(y)\hat{u}=0,
\nonumber
\end{eqnarray}
where 
\begin{equation}                               \label{eq:2.8}
V_1=(A_n')^2+i\frac{\partial A_n'}{\partial y_n}
+\sum_{j,k=1}^{n-1}\left(\hat{g}^{jk}A_j'A_k'+
i\frac{\partial}{\partial y_j}(\hat{g}^{jk}A_k')\right)
+\hat{V}(y).
\end{equation}
Note that
$
\hat{u}(y,0)=\hat{u}_{t}(y,0)=0,
\ \ \ \hat{u}|_{\Gamma\times(0,T)}=f(y',t)
$
since $\varphi(x',0)=x'.$

Denote by $\hat{\Lambda}$ the D-to-N operator corresponding
to $\hat{L}$.  We have
\begin{equation}                                   \label{eq:2.9}
\hat{\Lambda}f=
\left(\frac{\partial\hat{u}(y',y_n,t)}{\partial y_n}
+i\hat{A}_n(y',y_n)\hat{u}(y',y_n,t)\right)\left|_{\Gamma\times(0,T_0)}\right. .
\end{equation}
Note that 
\begin{equation}                                  \label{eq:2.10}
\Lambda f=\hat{\Lambda}f \ \ \ \mbox{on}\ \ \Gamma\times(0,T_0),
\end{equation}
where $\mbox{supp\ }f\subset\Gamma\times(0,T_0)$.
 It will be shown in Remark 2.2 that 
the D-to-N operator $\hat{\Lambda}$ on $\Gamma\times(0,T_0)$
determines the restriction of the metric tensor $\|\hat{g}^{jk}\|^{-1}$
to $\Gamma$ and, in particular,
\begin{equation}                                 \label{eq:2.11}
\hat{g}(y',0) \ \ \ \mbox{and}\ \ \ \ \ 
\frac{\partial\hat{g}}{\partial y_n}(y',0).
\end{equation}

Make the transformation $\hat{u}=(\hat{g}(y',y_n))^{-\frac{1}{4}}u'$.
Then $L'u'=0$ in $U_0\times(0,T_0)$ where $L'$ is similar to
(\ref{eq:2.7}) with $\hat{A}_j+A_j'$ replaced by $\hat{A}_j,\ 1\leq j\leq n$.
Denote by $\Lambda'$ the D-to-N operator corresponding to $L'$,  i.e.
\begin{equation}                            \label{eq:2.12}
\Lambda'f'=\left(\frac{\partial u'}{\partial y_n}
+i\hat{A}_nu'\right)\left|_{\Gamma\times(0,T)}\right.
\ \ \ \ \mbox{where}\ \ \  f'=u'\left|_{\Gamma\times(0,T)}\right.
\end{equation}
We have
\begin{equation}                            \label{eq:2.13}
\Lambda'f'=\left(\frac{\partial }{\partial y_n}
+i\hat{A}_n\right)
((\hat{g})^{\frac{1}{4}}\hat{u})
\left|_{\Gamma\times(0,T)}\right.
=\hat{g}^{\frac{1}{4}}(y',0)\hat{\Lambda} f 
+\frac{\hat{g}_{y_n}(y',0)}{4\hat{g}(y',0)}f'(y',0,t).
\end{equation}
It follows from (\ref{eq:2.11}) that $\hat{\Lambda}$ determines
$\Lambda'$ on $\Gamma\times (0,T_0)$.

Let $\hat{\psi}\in C^\infty(\overline{U_0}),\ \hat{\psi}(y',0)=0$.  
Then $c(y',y_n)=e^{i\hat{\psi}}\in G_0(\overline{U_0})$.
We say that 
$A_j^{(1)}(y)$ and $\hat{A}_j(y)$ are gauge equivalent 
in $\overline{U_0}$
if
\begin{equation}                                  \label{eq:2.14}
A_j^{(1)}=\hat{A}_j+\frac{\partial\hat{\psi}}{\partial y_j}
\ \ \ \mbox{in}\ \ \ \overline{U_0},\ \ 1\leq j\leq n.
\end{equation}

We shall choose $\hat{\psi}(y)$ such that $A_n^{(1)}=0$.  Denote
$u'=e^{i\hat{\psi}}\hat{u}_1$.  Then $\hat{u}_1$ satisfies the equation
\[
L_1\hat{u}_1=0,
\]
where 
\begin{eqnarray}                               \label{eq:2.15}
L_1\hat{u}_1=\frac{\partial^2\hat{u}_1}{\partial t^2}-
\frac{\partial^2\hat{u}_1}{\partial y_n^2}
\ \ \ \ \ \ \ \ \ \ \ \ \ \ \ \ \ \ \ \ \ \ \ \ \ \  \ \ \ \ \  
\\
 +
\sum_{j,k=1}^{n-1}\left(-i\frac{\partial}{\partial y_j}+A_j^{(1)}(y)\right)
\hat{g}^{jk} \left(-i\frac{\partial}{\partial y_k}+A_k^{(1)}(y)\right)\hat{u}_1
+V_1(y)\hat{u}_1=0,
\nonumber
\end{eqnarray}
i.e.
$L_1$ is obtained from $L'$ by replacing $\hat{A}_j$ by $A_j^{(1)}$.
Note that the D-to-N operator $\Lambda^{(1)}$ corresponding to $L_1$
has the form:
\begin{equation}                                 \label{eq:2.16}
\Lambda^{(1)}f=\frac{\partial \hat{u}_1}{\partial y_n}\left|_{y_n=0,\ 
0<t<T_0}\right.,
\end{equation}
since $A_n^{(1)}=0$.  Note that
\begin{equation}                                \label{eq:2.17}
\Lambda^{(1)}f=\Lambda'f,
\end{equation}
where $\Lambda'$ is the D-to-N corresponding to $L'$.
Denote by $L_1^*$ the formally adjoint operator to $L_1$.
Then $L_1^*$ has the same form as (\ref{eq:2.15})  with
$A_j^{(1)},V_1$ replaced by $\overline{A_j^{(1)}},\overline{V_1}$.  Let
$\Lambda_*^{(1)}$ be the D-to-N operator corresponding to $L_1^*$.
Note that $\Lambda^{(1)}$ determines $\Lambda_*^{(1)}$.
Indeed let $\Lambda^{(1)*}$ be the adjoint operator to $\Lambda^{(1)}$.
Then $\Lambda_*^{(1)}f$ is obtained from $\Lambda^{(1)*}$ by changing $t$
to $T-t$ and $f(y',t)$ to $f(y',T-t)$ (c.f. [KL1]).
Note that $L_1$ is self-adjoint if $L$ is self-adjoint.  In this case 
$\Lambda_*^{(1)}=\Lambda^{(1)}$. 

Denote $\Delta_{1s_0}=\Gamma\times(s_0,T]$, where $0\leq s_0<T$.  Let
$D(\Delta_{1s_0})$ be the forward domain of influence of $\Delta_{1s_0}$
in the half-space $y_n\geq 0$.  Denote 
$\Gamma^{(1)}=\{y':(y',y_n,t)\in D(\Delta_{10}),y_n=0,t=T\}$ and
$\Delta_{2s_0}=\Gamma^{(1)}\times (s_0,T]$.
Let
$D(\Delta_{2s_0})$ be the domain of influence of $\Delta_{2s_0}$ in 
$y_n\geq 0$.
Denote by $Y_{js_0}$ the intersection of $D(\Delta_{js_0})$ with the plane
$T-t-y_n=0$ and by $X_{js_0}$ the part of $D(\Delta_{js_0})$ below $Y_{js_0},
j=1,2.$  
Let $Z_{js_0}=\partial X_{js_0}\setminus(Y_{js_0}\cup\{y_n=0\}).$  Also denote 
$\gamma_{js_0}=\partial Y_{js_0}\cap\{y_n=0\},\ 1\leq j\leq 2.$
Note that $\gamma_{10}=\Gamma^{(1)}$.
Denote by $R_{s_0}$ the following subset of 
$Y_{2s_0}:\{(s,y'):s_0\leq s \leq T, y'\in \Gamma^{(1)}\}$,
$s=t-y_n,T-y_n-t=0$.
We shall choose $T$ such that $D(\Delta_{20})$
does not intersect $\partial\Omega \times[0,T]$ when $y_n>0$.
We also assume that the semigeodesic coordinates are defined in
$D(\Delta_{20})$ 
for $t\leq T$ and that $D(\Delta_{20})\cap\{y_n=0\}\subset \Gamma_0\times [0,T]$
for $t\leq T.$

The main result of this section is the following lemma: 
\begin{lemma}                              \label{lma:2.1}
Let $L_1$ be self-adjoint.  Then
the D-to-N operator $\Lambda^{(1)}$ on $\Delta_{20}$
determines all coefficients of $L_1$ (see (\ref{eq:2.15}) )
in $\D_{\frac{T}{2}}=\Gamma\times[0,\frac{T}{2}]$.
\end{lemma}

For the most part of the proof of Lemma \ref{lma:2.1} we will
not need $L_1$ to be self-adjoint.  We shall use the self-adjointness only 
in the proof of Lemma \ref{lma:2.4}.

Consider the identity
\begin{equation}                            \label{eq:2.18}
0=(L_1\hat{u}_1,\frac{\partial\hat{v}_1}{\partial t})
+(\frac{\partial \hat{u}_1}{\partial t},L_1^*\hat{v}_1),
\end{equation}
where
\begin{equation}                           \label{eq:2.19}
(\hat{u},\hat{v})=
\int_{X_{20}}\hat{u}(y,t)\overline{\hat{v}(y,t)}dydt,
\end{equation}
$L_1$ is the same as in (\ref{eq:2.15}), 
$\hat{u}_1(y',0,t)=f(y',t),\ 
\hat{v}_{1}(y',0,t)=g(y',t),\ \hat{u}=\hat{u}_t=
\hat{v}=\hat{v}_t=0$ for $t=0, y_n>0, L_1\hat{u}_1=0,\ L_1^*\hat{v}_1=0,
\ \mbox{supp\ }f(y',t)\subset \Delta_{20},\ 
\ \mbox{supp\ }g(y',t)\subset\Delta_{20}$.
We have:
\begin{eqnarray}                            \label{eq:2.20}
\int_{X_{20}}(\hat{u}_{1tt}\overline{v_{1t}} +\hat{u}_{1t}
\overline{\hat{v}_{1tt}})dydt 
= \int_{X_{20}}\frac{\partial}{\partial t}
(\hat{u}_{1t}\overline{\hat{v}_{1t}})dydt
\nonumber
\\
=\int_{Y_{20}}\hat{u}_{1t}\overline{\hat{v}_{1t}}dy
-\int_{Z_{20}}\hat{u}_{1t}\overline{\hat{v}_{1t}}dy
\end{eqnarray}
Since $\hat{u}_1$ and $\hat{v}_1$ have zero Cauchy data for
$t=0,\ y_n>0,$  we get that $\hat{u}_1,\ \hat{v}_1$ are
equal to zero on $Z_{20}$.
Integrating by parts in $y_n$ we obtain:
\begin{eqnarray}                               \label{eq:2.21}
-\int_{X_{20}}(\hat{u}_{1y_n^2}\overline{\hat{v}_{1t}}
+\hat{u}_{1t}\overline{\hat{v}_{1y_n^2}})dy'dy_ndt
=
\int_{X_{20}}(\hat{u}_{1y_n}\overline{\hat{v}_{1y_nt}}
+\hat{u}_{1y_nt}\overline{\hat{v}_{1y_n}})dy'dy_ndt
\nonumber
\\
\ \ \ \ \ \ \ \ \ -\int_{Y_{20}}(\hat{u}_{1y_n}\overline{\hat{v}_{1t}}
+\hat{u}_{1t}\overline{\hat{v}_{1y_n}})dy'dt
+\int_{\Delta_{20}}(\hat{u}_{1y_n}\overline{\hat{v}_{1t}}
+\hat{u}_{1t}\overline{\hat{v}_{1y_n}})dy'dt.
\end{eqnarray}
Note that 
\[
\hat{u}_{1y_n}\overline{\hat{v}_{1y_nt}}
+\hat{u}_{1y_nt}\overline{\hat{v}_{1y_n}}
=\frac{\partial}{\partial t}(\hat{u}_{1y_n}\overline{\hat{v}_{1y_n}}).
\]
Analogously, integrating by parts in $y'=(y_1,...,y_{n-1})$
we get
\begin{eqnarray}                               \label{eq:2.22}
\int_{X_{20}}\left[\sum_{j,k=1}^{n-1}\left(-i\frac{\partial}{\partial y_j}
+A_{j}^{(1)}\right)\hat{g}^{jk}\left(-i\frac{\partial}{\partial y_k}
+A_{k}^{(1)}\right)\hat{u}_1\overline{\hat{v}_{1t}}  
\right.
\\
+ 
\left. 
\hat{u}_{1t}\sum_{j,k=1}^{n-1}\left(i\frac{\partial}{\partial y_j}
+A_{j}^{(1)}\right)\hat{g}^{jk}\left(i\frac{\partial}{\partial y_k}
+A_{k}^{(1)}\right)\overline{\hat{v}_{1}}  \right] dydt
+\int_{X_{20}}V_1\hat{u}_{1t}\overline{v_1}dydt
+\int_{X_{20}}V_1\hat{u}_{1}\overline{v_{1t}}dydt
\nonumber
\\
=\int_{X_{20}}\frac{\partial}{\partial t}
\sum_{j,k=1}^{n-1}\hat{g}^{jk}(y)\left(-i\frac{\partial}{\partial y_k}
+A_{k}^{(1)}\right)\hat{u}_1\overline{\left(-i\frac{\partial}{\partial y_j}
+\overline{A_{j}^{(1)}}\right)\hat{v}_1}dydt 
+\int_{X_{20}}V_1\frac{\partial}{\partial t}
(\hat{u}_1\overline{\hat{v}_1})dydt.
\nonumber
\end{eqnarray}
Let
\begin{equation}                               \label{eq:2.23}
t-y_n=s,\ \ \ T-t-y_n=\tau,
\end{equation}
i.e. $y_n=\frac{T-\tau-s}{2},\ \ t=\frac{T-\tau+s}{2}.$
Note that
\[
u_s=\frac{1}{2}(u_t-u_{y_n}),\ \ u_\tau =-\frac{1}{2}(u_t+u_{y_n})
\]
and that  $\tau=0$ on $Y_{20}$.  Combining
(\ref{eq:2.20}), (\ref{eq:2.21}), (\ref{eq:2.22}) we get
\begin{equation}                             \label{eq:2.24}
0=(L_1\hat{u}_1,\hat{v}_{1t})+(\hat{u}_{1t},L_1^*\hat{v}_1)
=Q(\hat{u}_1,\hat{v}_1)-
\Lambda_0(f,g),
\end{equation}
where
\begin{eqnarray}                             \label{eq:2.25}
\ \ \ \ \ \ \ \ \ \ \ \ \ \ \ \ \ \ \ \ \ \ \ \ \ \ \ \ \ \ \   
Q(\hat{u}_1,\hat{v}_1)\ \ =
\ \ \ \ \ \ \ \ \ \ \ \ \ \ \ \ \ \ \ \ \ \ \ \ \ \ \ \ \ \ \ \ 
\ \ \ \ \ \ \ \ \ \ \ \ \ \ \ \ \ \ \ \ \ \ \  
\\
\int_{Y_{20}}\frac{1}{2}\left[4\hat{u}_{1s}\overline{\hat{v}_{1s}}+
\sum_{j,k=1}^{n-1}\hat{g}^{jk}\left(-i\frac{\partial}{\partial y_j}
+A_{j}^{(1)}\right)\hat{u}_1\overline{\left(-i\frac{\partial}{\partial y_k}
+\overline{A_{k}^{(1)}}\right)\hat{v}_1}
+
V_1(y)\hat{u}_1\overline{\hat{v}_1}\right]dy'ds ,
\nonumber
\end{eqnarray}

\begin{equation}                          \label{eq:2.26}
\Lambda_0(f,g)=-\int_{\Delta_{20}}\hat{u}_{1y_n}\overline{\hat{v}_{1t}}+
\hat{u}_{1t}\overline{\hat{v}_{1y_n}}dy'dt.
\end{equation}
Note that 
$u_{1t}|_{\Delta_{20}}=f_t,\ \hat{v}_{1t}|_{\Delta_{20}}=g_{1t},
\ \hat{u}_{1y_n}|_{\Delta_{20}}=\Lambda^{(1)}f,
\ v_{1y_n}|_{\Delta_{20}}=\Lambda_*^{(1)}g$,
where $\Lambda^{(1)}$ is the D-to-N operator on $
\Delta_{20}$.

Therefore $\Lambda_0(f,g)$ is known for all $f,g$ in
$H_0^1(\Delta_{20})$ if we know $\Lambda^{(1)}$.  Here
$H_0^1(\Delta_{js_0})$ is the closure in 
$H^1(\Delta_{js_0})$ of 
smooth functions equal to zero on 
$\partial\Delta_{js_0}\setminus\{t=T\},\ 0\leq s_0 < T,\ j=1,2.$

For the convenience we shall often use the notation $u^f,\ 
v^g$ 
where $L_1u^f=0$ and $L_1^*v^g=0$ in $X_{20}$
to emphasize the dependence of $u^f$ and
$v^g$ on $f$ and $g$ respectively.
Note that $Q(\hat{u}_1,\hat{v}_1)$ is a bounded nonsymmetric bilinear
form on $H_0^1(Y_{20})$ where $H_0^1(Y_{js_0})$ consists of functions in
$H^1(Y_{js_0})$ equal to zero on $\partial Y_{js_0}\setminus \{t=T\},
\ 0\leq s_0<T,\ 1\leq j\leq 2$.

Note that for $T$ small we have for any $\hat{u}_1\in H_0^1(Y_{20})$
\begin{equation}                             \label{eq:2.27}
\Re Q(\hat{u}_1,\hat{u}_1)\geq C\|\hat{u}_1\|_1^2,
\end{equation}
where $\|\ \ \|_1$ is the norm in $H^1(Y_{20})$.
Also we have
\begin{equation}                          \label{eq:2.28}
|Q(\hat{u}_1,\hat{v}_1)|=|\Lambda_0(f,g)|
\leq C\|f\|_1\|g\|_1,
\end{equation}
where $\|\ \ \|_1$ is the norm in $H^1(\Delta_{20})$.

Note that on the plane $\tau=0$ we have 
$\Gamma\times[s_0,T]\subset \overline{Y_{1s_0}}
\subset R_{s_0}\subset \overline{Y_{2s_0}}$
where $R_{s_0}=
\Gamma^{(1)}\times[s_0,T]$.
\begin{lemma}                           \label{lma:2.2}
For any $w\in H_0^1(R_{s_0})\subset H_0^1(Y_{2s_0})$ there exists
a sequence $\{u^{f_n}\},f_n\in H_0^1(\Delta_{2s_0})$  that converges to $w$
in $H_0^1(Y_{2s_0})$.
\end{lemma}
We shall postpone the proof of Lemma \ref{lma:2.2} until
the section 3.

Let the span of $u_j,\ j\geq 1,$ be dense in
$H_0^1(Y_{2s_0})$.
\begin{lemma}                              \label{lma:2.3}
(c.f. [GLT])\ 
For any $u\in H_0^1(Y_{20})$ there exists $u_0\in H_0^1(Y_{2s_0})$
such that $Q(u_0,v')=Q(u,v'),\ \forall v'\in H_0^1(Y_{2s_0})$.
\end{lemma} 
{\bf Proof}:
Take $u^{(N)}=\sum_{j=1}^N c_{jN}u_j$ and determine $c_{jN}$ from
the linear system
\begin{equation}                            \label{eq:2.29}
Q(u^{(N)},u_k)=Q(u,u_k),\ \ 1\leq k<N,
\end{equation}
i.e.
\begin{equation}                            \label{eq:2.30} 
\sum_{j=1}^Nc_{jN}Q(u_j,u_k)=Q(u,u_k),\ \ 1\leq k\leq N.
\end{equation}
  Multiplying (\ref{eq:2.29})
by $\overline{c_{kN}}$ and adding,  we get that
\begin{equation}                            \label{eq:2.31}
Q(u^{(N)},u^{(N)})=Q(u,u^{(N)}).
\end{equation}
Therefore
\[
C_1\|u^{(N)}\|_1^2\leq \Re Q(u^{(N)},u^{(N)})\leq C\|u\|_1\|u^{(N)}\|_1.
\]
Therefore system (\ref{eq:2.29}) has a unique solution and 
$\|u^{(N)}\|_1\leq C_2\|u\|_1$ for all $N$.
Since the sequence $\|u^{(N)}\|_1$ is bounded it is weakly compact in 
$H_0^1(Y_{2s_0})$.  Therefore there exists a subsequence $u^{(N_k)}$ such
that $u^{(N_k)}\rightarrow u_0\in H_0^1(Y_{2t_0})$ weakly when 
$N_k\rightarrow\infty$.
Passing to the limit in (\ref{eq:2.29})
we get
\begin{equation}                            \label{eq:2.32}
Q(u_0,u_k)=Q(u,u_k),\ \ \forall k.
\end{equation}
We used that $Q$ is continuous bilinear form in $H_0^1(Y_{2s_0})$.  
Since the span of $\{u_k\}$ is dense we get 
$Q(u_0,v')=Q(u,v'),\ \forall v'\in H_0^1(Y_{2s_0})$.
Note that any subsequence of $\{u^{(N)}\}$ has a subsequence that
converges to $u_0$ since (\ref{eq:2.27})  implies that there is a unique
$u_0$ satisfying (\ref{eq:2.32}).  Therefore the whole sequence 
$u^{(N)}$ converges weakly to $u_0$.  Note that $u^{(N)}$ converges 
also strongly in $H_0^1(Y_{2s_0})$ to $u_0$ since
\[
Q(u_0-u^{(N)},u_0-u^{(N)})=Q(u_0,u_0)-Q(u_0,u^{(N)})-
Q(u^{(N)},u_0)+Q(u^{(N)},u^{(N)}).
\]
Since $Q(u^{(N)},u^{(N)})=Q(u,u^{(N)})$ and $Q(u_0,u_0)=Q(u,u_0)$
we get that
$Q(u_0-u^{(N)},u_0-u^{(N)})\rightarrow 0$ when $N\rightarrow \infty$.
\qed

Note that Lemma \ref{lma:2.3} remains true when $Q(u,v')$ is replaced
by any linear continuous functional $\Phi(v')$ on $H_0^1(Y_{2s_0})$.

Applying the Green's formula in $X_{20}$ we get,  similarily to
(\ref{eq:2.24}):
\begin{eqnarray}                         \label{eq:2.33}
0=\int_{X_{20}}(L_1u\overline{v}-u\overline{L_1^*}v)dydt
\\
=\int_{Y_{20}}(u_s\overline{v}-u\overline{v_s})dsdy'
-\int_{\Delta_{20}}(u_{y_n}\overline{v}-u\overline{v_{y_n}})dy'dt.
\nonumber
\end{eqnarray}
We used in (\ref{eq:2.33})
that $u=v=0$ on $Z_{20}$.  Integrating by parts we get
\begin{equation}                          \label{eq:2.34}
-\int_{Y_{20}}u\overline{v_s}dsdy'=
-\int_{\gamma_{20}}u(y',0,T)\overline{v(y',0,T)}dy' +
\int_{Y_{20}}u_s\overline{v}dsdy'.
\end{equation}
Let
\begin{equation}                          \label{eq:2.35}
A(u,v)=2\int_{Y_{20}}u_s\overline{v}ds dy'.
\end{equation}
Since $u(y',0,T)=f(y',T)$ and $v(y',0,T)=g(y',T)$ we get
from (\ref{eq:2.33}) and (\ref{eq:2.34}) that $A(u,v)$ is 
determined by the boundary data.

Fix $f,g$ in $H_0^1(\Delta_{10})$ and smooth.  Let 
$s_0\in [0,T)$ be arbitrary.  We shall show that there exists a unique
$u_0\in H_0^1(R_{s_0})$ such that
\begin{equation}                            \label{eq:2.36}
A(u_0,v')=A(u^f,v'),\ \ \forall v'\in H_0^1(Y_{2s_0}).
\end{equation}
Note that $\overline{Y_{10}}\cap\{s\geq s_0\}\subset R_{s_0}$. 
Denote by $w(s,y')$ a function equal to $u^f|_{s=s_0}$  
in 
$R_{s_0}$ and
 and equal to zero in $Y_{2s_0}\setminus R_{s_0}$.  Then 
$\frac{\partial w}{\partial s}=0$ in $Y_{2s_0},\ u_0=u^f-w\in H_0^{1}(R_{s_0})$ 
and $\frac{\partial u^f}{\partial s}=\frac{\partial u_0}{\partial s}$
in $Y_{2s_0}$.  Therefore (\ref{eq:2.36}) holds.

We shall show that $A(u_0,v^g)$ is uniquely determined by the D-to-N operator.

\begin{lemma}                             \label{lma:2.4}
Let $L_1^{(i)},\ i=1,2$, be two formally self-adjoint operators 
in $X_{20}^{(i)},\ i=1,2$,  
such that the corresponding D-to-N operators $\Lambda_1^{(i)},\ i=1,2$,
are equal on $\Delta_{2s_0}$. Then
\begin{equation}
C_1\|u_1^f\|_{1,Y_{2s_0}^{(1)}}\leq \|u_2^f\|_{1,Y_{2s_0}^{(2)}}\leq 
C_2\|u_1^f\|_{1,Y_{2s_0}^{(1)}}
\nonumber
\end{equation}
for all $f\in H_0^1(\Delta_{2s_0})$.
\end{lemma}

Here $L_1^{(i)}u_i^f=0$ in $X_{20}^{(i)}$.

We shall show first that $\Delta_{2s_0}^{(1)}=\Delta_{2s_0}^{(2)}$.
Denote by $\Delta^{(i)}(0,T)$ the intersection of the domain
of influence $D(\Delta_{10})$ of $L_1^{(i)}$ with $y_n=0,i=1,2.$
Note that $\Delta^{(i)}(0,T)
=(\overline{\cup_{f}\mbox{supp}\ u_i^f})\cap \{y_n=0\}$ 
where the union is taken over all $f\in H_0^1(\Delta_{10})$.
Denote $\tilde{\Delta}^{(i)}=(\Gamma\times[0,T])\cup G_i$  where $G_i$ is
the closure of the union over all $f\in H_0^1(\Delta_{10})$ of 
$\mbox{supp} \ \Lambda^{(i)}f,i=1,2$.
It is clear that 
$\tilde{\Delta}^{(i)}\subset \Delta^{(i)}(0,T),i=1,2$.  The inclusion 
$ \Delta^{(i)}(0,T)\subset\tilde{\Delta}^{(i)}$
follows from the local uniqueness of the Cauchy problem  (see [T]).  Since 
$\Lambda^{(1)}=\Lambda^{(2)}$  we have  $\tilde{\Delta}^{(1)}=
\tilde{\Delta}^{(2)}$.  Therefore  $\Delta^{(1)}(0,T)=\Delta^{(2)}(0,T)$.
In particular,  $\Gamma_1^{(1)}=\Gamma_2^{(1)}$ and therefore 
$\Delta_{2s_0}^{(1)}=\Delta_{2s_0}^{(2)}$.  Now we shall proceed with the 
proof of Lemma \ref{lma:2.4}.

Consider the identity (\ref{eq:2.24}) for $L_1^{(1)}$ and $L_1^{(2)}$ 
respectively assuming thet  
$(L_1^{(i)})^*=L_1^{(i)}$ and $u_i^f=v_i^f,\ i=1,2$.  Since
$\Lambda_0^{(1)}(f,f)=\Lambda_0^{(2)}(f,f)$ we have that
$Q^{(1)}(u_1^f,u_1^f)=Q^{(2)}(u_2^f,u_2^f)$  where the quadratic form
$Q^{(i)}$ corresponds to $L_1^{(i)}$.  Since $T$ is small
$Q^{(i)}(u_i^f,u_i^f)$ is equivalent to $\|u_i^f\|_{1,Y_{2s_0}^{(i)}}^2$.
\qed

We shall show that Lemma \ref{lma:2.4} implies that 
\begin{equation}                      \label{eq:2.37}
A^{(1)}(u_0^{(1)},v_1^g)=A^{(2)}(u_0^{(2)},v_2^g)
\end{equation}
  where $A^{(i)}$ and
$u_0^{(i)}$ correspond to $L_1^{(i)},\ i=1,2$.  It follows from
(\ref{eq:2.33}), (\ref{eq:2.34})  that $A^{(1)}(u_1^{f'},v_1^{g'})=
A^{(2)}(u_2^{f'},v_2^{g'})$ for any $f'\in H_0^1(\Delta_{2s_0}),\ 
g'\in H_0^1(\Delta_{20})$.

By Lemma \ref{lma:2.2}  there exists a sequence 
$f_n\in H_0^1(\Delta_{2s_0})$
such that $\|u_0^{(1)}-u_1^{f_n}\|_{1,Y_{2s_0}^{(1)}}\rightarrow 0$
when $n\rightarrow\infty$.
By Lemma \ref{lma:2.4} the sequence $\{u_2^{f_n}\}$ is also convergent
in $H_0^1(Y_{2s_0}^{(2)})$ to some 
$w^{(2)}\in H_0^1(Y_{2s_0}^{(2)})$.  Since $A^{(i)}(u_i^{f'},v_i^{g'})$
are continuous functionals on $H_0^1(Y_{2s_0}^{(i)}),\ i=1,2,$
we get passing to the limit that
\begin{equation}                     \label{eq:2.38}
A^{(1)}(u_0^{(1)},v_1^{g'})=A^{(2)}(w^{(2)},v_2^{g'}).
\end{equation}
In (\ref{eq:2.38})  $g'$ is arbitrary  in $H_0^1(\Delta_{20})$.
Take $g'\in H_0^1(\Delta_{2s_0})$.
Then (\ref{eq:2.36}) implies that
\begin{equation}                     \label{eq:2.39}
A^{(1)}(u_0^{(1)},v_1^{g'})=A^{(2)}(u_0^{(2)},v_2^{g'}).
\end{equation}
 Compairing (\ref{eq:2.38}) and (\ref{eq:2.39}) 
we get $A^{(2)}(u_0^{(2)},v_2^{g'})=
A^{(2)}(w^{(2)},v_2^{g'})$.  By Lemma \ref{lma:2.2} 
$\{v_2^{g'}\}, g'\in H_0^1(\Delta_{2s_0})$ are dense in 
$H_0^1(R_{s_0})$.  
Therefore $u_0^{(2)}=w^{(2)}$ in $R_{s_0}$.  We assume that 
$v^g\in H_0^1(Y_{10})$.
Since $Y_{10}\cap \{s\geq s_0\}\subset R_{s_0}$  we get that 
$A^{(2)}(w^{(2)},v^g)=A^{(2)}(u_0^{(2)},v^g)$. 
Therefore (\ref{eq:2.37}) holds,  i.e. $A(u_0,v^g)$ is uniquely determined
by the D-to-N operator.

{\bf Remark 2.1}
In this remark we shall show that in the self-adjoint case
$A(u_0,v^g)$ can be recovered constructively from the boundary data.
For any $\e>0$
$Q_\e(u^f.v^g)=\e Q(u^f,v^g)+A(u^f,v^g)$ is determined 
by the D-to-N operator.  
Let $g_j,\ j\geq 1$ be a dense set in $H_0^1(\Delta_{2s_0})$.
Denote by $\overline{H}$ the closure of the span of 
$\{v^{g_j}\}$ in $H_0^1(Y_{2s_0})$.
Since $\Re A(v',v')\geq 0$  (see (\ref{eq:2.35}) ) we get
that $\Re Q_\e(v',v')\geq C\e \|v'\|_1^2$ for any $v'\in H_0^1(Y_{2s_0})$.
It follows from (\ref{eq:2.30}) with $Q$ replaced by $Q_\e$, 
$u_j$ replacedby $v^{g_j}$
and $c_{jN}$
replaced by $c_{jN\e}$ that $c_{jN\e}$ are determined by the D-to-N
operator.  Repeating the proof of Lemma \ref{lma:2.3} we get that 
$u_\e^{(N)}=\sum_{j=1}^N c_{jN\e}v^{g_j}$ converges in $H_0^1(Y_{2s_0})$ to 
$u_\e\in \overline{H}$ such that $Q_\e(u_\e,v')=Q_\e(u^f,v'),
\ \forall v'\in \overline{H}$.
Since $Q_\e(u_\e^{(N)},v^g)$ is determined by the D-to-N operator 
the limit $Q_\e(u_\e,v^g)=\lim_{N\rightarrow\infty}Q_\e(u_\e^{(N)},v^g)$ 
is also determined
by the D-to-N operator.  
Here $v^g\in H_0^1(Y_{10})$.
Denote  $w_\e=u_\e-u_0$,  where $u_0$ is the same
as in (\ref{eq:2.36}),  $u_0\in H_0^1(R_{s_0})$.  
We have,  using (\ref{eq:2.36}),  that
\begin{equation}
\e Q(w_\e,v')+A(w_\e,v')=\e Q(u^f-u_0,v').
\nonumber
\end{equation} 
Take $v'=w_\e$.  Then
\begin{equation}
C\e\|w_\e\|_1^2\leq C\e \|u^f-u_0\|_1\|w_\e\|_1.
\nonumber
\end{equation} 
Therefore $\|w_\e\|_1\leq C\|u^f-u_0\|_1$ for all $\e$.  As in the proof 
of Lemma \ref{lma:2.3} we get that a sequence $w_{\e_k}$ converges
weakly in $\overline{H}$ to some $w\in\overline{H}$ and $(w_s,v')=0\ \ 
\forall v'\in\overline{H}$.  Since
$\overline{H}\supset H_0^1(R_{s_0})$  we get that $w=0$ in $R_{s_0}$.  
Since $R_{s_0}\supset Y_{10}\cap\{s\geq s_0\}$ we have that
$(w_s,v^g)=0$.
Therefore $Q_\e(u_\e,v^g)\rightarrow A(u_0,v^g)$.  Since this is true 
for any sequence $\e_k$  
we have $A(u_0,v^g)=\lim_{\e\rightarrow 0} Q_\e(u_\e,v^g)$,
i.e. $A(u_0,v^g)$ is determined by the D-to-N operator on $\Delta_{20}$.
\qed

Denote
\[
A_1(u^f,v^g)=A(u^f,v^g)-A(u_0,v^g).
\]
Then $A_1(u^f,v^g)$ is also determined by the D-to-N operator and

\begin{equation}                      \label{eq:2.40}
A_1(u^f,v^g)=2\int_{Y_{10}\cap\{s\leq s_0\}}
\frac{\partial u^f}{\partial s}\overline{v^g}dsdy',
\end{equation}
since $u_s^f-u_{0s}=0$ when $s\geq s_0,\ u_0=0$ when $s\leq s_0$.

Now we shall construct a geometric optics solution $\hat{u}(y)$ of
 $L_1\hat{u}=0$
such that $\hat{u}(y,0)=\hat{u}_t(y,0)=0$ for $y_n>0$ and substitute it
in (\ref{eq:2.40}) to recover $v^g$.  We are looking for $\hat{u}$
in the form:
\begin{equation}                      \label{eq:2.41}
\hat{u}=u_N+u^{(N+1)},
\end{equation}
where 
\[
u_N=e^{ik(s-s_0)}\sum_{p=0}^N\frac{1}{(ik)^p}a_p(s,\tau,y').
\]
Substituting $\hat{u}$ in 
$L_1\hat{u}=(-4\frac{\partial^2}{\partial s\partial\tau}+L_1')\hat{u}=0$ we get
the transport equations for $a_p$:
\[
\frac{\partial a_0}{\partial \tau}=0,
\ \ 4\frac{\partial a_p}{\partial \tau}=L_1a_{p-1},\ \ p\geq 1.
\]
Choose $a_0(s,y')=\chi_1(s)\chi_2(y')$,  where 
$\chi_1(s)\in C_0^\infty(\R^1),\ \chi_1(s)=1$ for $|s-s_0|<\delta,\ 
\chi_1(s)=0$ for $|s-s_0|>2\delta,
\ \chi_2(y')=\frac{1}{\e^{n-1}}\chi_0(\frac{y'-y_0'}{\e}),\  
\chi_0(y')\in C_0^\infty(\R^{n-1}),\ \chi_0(y')=0$ for $|y'|>\delta,\ 
\int_{\R^{n-1}}\chi_0(y')dy'=1,\ \delta$ is small,
$y_0'\in \Gamma$.

We define $a_p=
\frac{1}{4}
\int_{T-s}^\tau(L_1a_{p-1})d\tau',\ 1\leq p\leq N,$ and
$u^{(N+1)}$ as the solution of 
\begin{eqnarray}
L_1u^{(N+1)}=
-\frac{1}{4^N(ik)^N}(L_1a_N)e^{ik(s-s_0)},
\nonumber
\\ 
\ \ \ \ \ \ \ \ u^{(N+1)}=  u_t^{(N+1)}=0 \ 
\mbox{when}\ \  t=0,
\nonumber
\\
\ u^{(N+1)}=0\ \ \ 
\mbox{when}\ \  y_n=0.
\nonumber
\end{eqnarray}
Note that $\mbox{supp}\ u_N$ is contained in a small neighborhood
of the line $\{s=s_0,y'=y_0'\}$.  Therefore $\mbox{supp}\ (u_N +u^{(N+1)})\cap \{ y_n=0\}
\subset \Delta_{10}$.

Substitute (\ref{eq:2.41}) into (\ref{eq:2.40}).
Note that the principal term in $k$ has the form
\[
ik\int_{s<s_0} e^{ik(s-s_0)}\chi_1(s)\chi_2(y')
\overline{v^g(y',s)}dy'ds,
\]
where $v^g(y',s)$ is $v^g(y',s,\tau)$ for $\tau=0$.
Integrating by parts in $s$ and taking the limit when 
$k\rightarrow\infty$ we get that the boundary data determine

\begin{equation}                             \label{eq:2.42}
\int_{\R^{n-1}}\chi_2(y')\overline{v^g(y',s_0)}dy'.
\end{equation}

Taking limit in (\ref{eq:2.42}) when $\e\rightarrow 0$ we can recover
$v^g(y_0',s_0),\ y_0'\in \Gamma, \ 0<s_0<T,\ \tau=0$.  
Changing $T$ to $T-\tau',\ 0<\tau'<T$, we can analogously recover
$v^g(y',y_n,t)$ 
for $y'\in \Gamma,\ 0\leq y_n\leq \frac{T}{2},\ y_n
\leq t\leq T-y_n$.
In particular we determine  $\hat{v}^g(y,t)$ 
and its time derivatives
for $t=\frac{T}{2}$
and 
$y\in\D_{\frac{T}{2}}$ where $\D_{\frac{T}{2}}=\Gamma\times[0,\frac{T}{2}]$.

It is known (see [B1]) that $\hat{v}^g(y,\frac{T}{2}),\ 
g\in C_0^\infty(\Delta_{10})$ are dense in $H^m(\D_{\frac{T}{2}})$
for any $m\geq 0$.
Since $L_1v^g=0$ we have
\begin{eqnarray}                                \label{eq:2.43}
\frac{\partial^2}{\partial t^2}\hat{v}^g(y,t)-
\frac{\partial^2\hat{v}^g}{\partial y_n^2}
-\sum_{j,k=1}^{n-1}\hat{g}^{jk}(y)\frac{\partial^2\hat{v}^g}
{\partial y_j\partial y_k}
\\
-\sum_{j=1}^{n-1}B_j\frac{\partial^2\hat{v}^g}{\partial y_j}
+C(y)\hat{v}^g=0,
\nonumber
\end{eqnarray}
where $B_j,\ 1\leq  j\leq n-1,\ C$ depend on
$\hat{g}^{jk},A_j^{(1)},V_1$ (see (2.15)).

Fix any $y=y_0\in \D_{\frac{T}{2}}$ and let $t=\frac{T}{2}$.  We can
consider (\ref{eq:2.43}) as a linear system with unknowns
$\hat{g}^{jk}(y_0),\ B_j(y_0),\ C(y_0),
\ g\in C_0^\infty(\Delta_{10})$
is arbitrary.  If the rank of (\ref{eq:2.43}) is not maximal
then there exist constants $\alpha_{jk},\alpha_j,\alpha_0,
\alpha_{jk}=\alpha_{kj}$, not
all equal to zero, such that
\begin{equation}                              \label{eq:2.44}
\sum_{j,k=1}^{n-1}\alpha_{jk}
\frac{\partial^2\hat{v}^g(y_0,\frac{T}{2})}
{\partial y_j\partial y_k}
+\sum_{j=1}^{n-1}\alpha_j\frac{\partial v^g(y_0,\frac{T}{2})}{\partial y_j}
+\alpha_0 v^g(y_0,\frac{T}{2})=0
\end{equation}
for all $\hat{v}^g(y,t),\ g\in C_0^\infty(\Delta_{10})$.  Then
we have a contradiction since $\{\hat{v}^g\}$ are  
dense in $H_0^N(\D_{\frac{T}{2}}),\ 
N\geq \frac{n}{2}+2$.  Therefore the system (\ref{eq:2.43})
has the maximal rank and $\hat{g}^{jk}(y_0),B_j(y_0),C(y_0)$ are uniquely
determined by (\ref{eq:2.43}).
In particular,  we recover  $\|\hat{g}^{jk}(y_0)\|$
in $\D_{\frac{T}{2}}$.  Knowing the metric and $B_j(y_0)$ we
recover $A_j^{(1)}(y_0),\ 1\leq j\leq n-1$.  Finally knowing 
$\hat{g}^{jk},A_j^{(1)}$ and $C(y_0)$ in $\D_{\frac{T}{2}}$ we can recover 
$\hat{V}(y_0)$.  Lemma \ref{lma:2.1}  is proven.
\qed

{\bf Remark 2.2}
We shall show that the D-to-N operator on $\Delta_{20}$ determines 
the metric tensor and (\ref{eq:2.11}) on $\Delta_{20}$.  Consider 
the cotangent space $T_0^*$ of $\Delta_{20}$.  Let $\xi_0$
be the dual variable to the time variable $t$.  Therefore the points
of $T_0^*$ have the form $(y',t,\xi',\xi_0)$.  The region 
$T_E^*\subset T_0^*$
 where $|\xi_0|<\e|\xi'|, \e$ is small, is a part of the "elliptic"
region of $T_0^*$ (see, for example, [H] or [E2])
where microlocally $L$ behaves as an elliptic operator.  Therefore
in the region $T_E^*$ we can find the parametrix for the D-to-N
operator $\Lambda$ as in the elliptic case (see, for example,  [LU] or 
[E5], pp 54-55).  
Therefore  we can recover the full symbol of $\Lambda$
in $T_E^*$,
in particular, 
   the principal symbol of 
$\Lambda$ in semigeodesic coordinates:
\begin{equation}
i
\left(
\sum_{j,k=1}^{n-1}g^{jk}(y',0)\xi_j\xi_k
-\xi_0^2
\right)^{\frac{1}{2}}
\nonumber
\end{equation}
It follows from [LU], formula (1.8), (see also [E5]),  that one can recover $\frac{\hat{g}_{y_n}(y',0)}{\hat{g}(y',0)}$.
\ 
\qed

\section{The conclusion of the proof of Lemma \ref{lma:2.1}.}
\label{section 3}
\init

We shall start with the proof of Lemma \ref{lma:2.2}.
We shall show first that 
it is enough to prove that the set $\{u^f\}$ where 
$f\in C_0^\infty(\Delta_{2s_0})$ is dense in 
$\stackrel{\circ}{H^1}(R_{s_0})$
where 
$\stackrel{\circ}{H^1}(Y_{2s_0})$
is the closure of $C_0^\infty(Y_{2s_0})$ in the $H^1(Y_{2s_0})$ norm
and $\stackrel{\circ}{H^1}(R_{s_0})$ is defined analogously.
Suppose there exists $v\in H_0^1(R_{s_0})$ that does not belong to
the closure $\overline{H}$ of $\{v^f\},g\in H_0^1(\Delta_{2s_0})$.
Then there exists $w\in H_0^1(Y_{2s_0})$  such that $(w,u^f)_1=0,\ 
\forall f\in H_0^1(\Delta_{2s_0})$, and
$(w,v)_1=1$,
 where $(v,w)_1$ is the inner product
in $H^1(Y_{2s_0})$.  

Take any $u_1\in C_0^\infty(R_{s_0})$.
Since by the assumption $\{u^f\},\ f\in C_0^\infty(\Delta_{2s_0})$
are dense in $\stackrel{\circ}{H^1}(R_{s_0})$ we get that
$(w,u_1)=0$.
Integrating by parts we get 
$\int_{R_{s_0}}(-\Delta w+w)\overline{u_1}dyds=0$.
Since $u_1\in C_0^\infty(R_{s_0})$ is arbitrary we get
$-\Delta w+ w=0$ in $R_{s_0}$.
Fix $y_1\in \Gamma^{(1)}$and consider $f\in H_0^1(\Delta_{2s_0})$ 
having support  in $\{|y'-y_1|<\e\}\times (T-\e,T]$,  where
$\e>0$ is small.  Then $\mbox{supp}\ u^f\subset R_{s_0}$. 
  Again integrating by parts we get
\[
0=(w,u^f)_1=\int_{R_{s_0}}(-\Delta w +w)\overline{u^f}dy'ds
\]
\[
+\int_{\Gamma^{(1)}}\frac{\partial w}{\partial s}\overline{u^f}(y',T)dy'
=\int_{\Gamma^{(1)}}\frac{\partial w}{\partial s}\overline{u^f}(y',T)dy'.
\]
Since $u^f=f(y',T)$ when $t=T$ and $f(y',T)$ can be chosen arbitrary
near $y_1$ we get that $\frac{\partial w}{\partial s}=0$  for $t=T,
|y'-y_1| < \e$.  Analogously we can prove that 
$\frac{\partial w}{\partial s}=0$ on $\Gamma^{(1)}\times\{t=T\}.$ 
  Then $(w,v)_1=\int_{R_{s_0}}(-\Delta w+w)\overline{v}dsdy'+
\int_{\Gamma^{(1)}}\frac{\partial w}{\partial s}\overline{v}dy'=0$
and this contradicts the assumption that $(w,v)_1=1$.   
\qed

To complete the proof of Lemma \ref{lma:2.2} we need two more lemmas.

Denote by $\Delta_1$ a domain in $\R^{n+1}$ bounded by three planes:
$\Gamma_2=\{\tau=T-t-y_n=0,\ 0\leq y_n\leq \frac{T}{2}\},
\ \Gamma_3=\{s=t-y_n=0,\ \frac{T}{2}\leq y_n\leq T\},
\ \Gamma_4=\{t=T,\ 0\leq y_n\leq T\}$.
Let $\mathcal{H}=\stackrel{\circ}{H^1}(\Gamma_4)\times L_2(\Gamma_4)$
and let $\mathcal{H}_1$ be the space of pairs $\{\varphi,\psi\},
\ \varphi\in H^1(\Gamma_2),\
\psi\in H^1(\Gamma_3),\ \varphi =0$   when $t=T,\ \psi=0$  when
$y_n=T,\ \varphi=\psi$ when $t=\frac{T}{2}$,  with the norm 
$\left(\|\varphi\|_{1,\Gamma_2}^2 + \|\psi\|_{1,\Gamma_3}^2\right)^{\frac{1}{2}}$.
We shall consider  functions with compact \mbox{supp}ort
in $y'\in \R^{n-1}$.

\begin{lemma}                                      \label{lma:3.1}
For any $\{v_0,v_1\}\in\mathcal{H}_1$
 there exist $\{w_0,w_1\}\in \mathcal{H}$
 and $u\in H^1(\Delta_1)$ such that 
$L_1u=0$ in $\Delta_1,\ u|_{\Gamma_2}=v_0,\ u|_{\Gamma_3}=v_1,
\ u|_{\Gamma_4}=w_0,\ \frac{\partial u}{\partial t}\left|_{\Gamma_4}\right.
=w_1$.
\end{lemma}

{\bf Proof:} For any smooth $u,v$ with compact \mbox{supp}orts in $y'$ such that 
$L_1u=0,L_1^*v=0$ in $\Delta_1$ we have
\begin{equation}                               \label{eq:3.1}
0=(u_t,L_1^*v)+(L_1u,v_t)=E(u,v)-Q(u,v)-Q_1(u,v),
\end{equation}
where $Q(u,v)$ is similar to (\ref{eq:2.25}),
$Q_1(u,v)$ is a bilinear form on $\Gamma_3$ of the form 
(\ref{eq:2.25}) with $u_s,\overline{v_s}$ replaced by
$u_\tau,\overline{v}_\tau$, 
 and 
\begin{eqnarray}
E(u,v)\ \
\ \ \ \ \ \ \ \ \ \ \ \ \ \ \ \ \ \ \ \ \ \ \ \ \ \ \ \ \ \ \ \
\ \ \ \ \ \ \ \ \ \ \ \ \ \ \ 
\nonumber
\\
=\int_{\Gamma_4}\left[ u_t\overline{v_t}+u_{y_n}\overline{v_{y_n}}
+\ \   \sum_{j,k=1}^{n-1}\hat{g}^{jk}(-i\frac{\partial}{\partial x_j}+
A_j^{(1)})
u
\overline{(-i\frac{\partial}{\partial x_k}+\overline{A_k^{(1)}})v}+
\hat{V}_1u\overline{v}\right]dy.
\nonumber
\end{eqnarray}

Consider $0=(L_1u,u_t)+(u_t,L_1u)$ in $\Delta_1$.  Integrate by parts
as in (\ref{eq:3.1}).  When $L_1$ is not self-adjoint we get,  in addition to
$Q(u,u), \ Q_1(u,u)$ and $E(u,u)$, an integral $E_1(u,u)$ over
$\Delta_1$ that satisfies
the following estimate:
\begin{equation}                                 \label{eq:3.2}
|E_1(u,u)|\leq C\int_{\Delta_1}(|u_t|^2+\sum_{j=1}^n|u_{y_j}|^2+|u|^2)dydt.
\end{equation}
Denote by $\Delta_{1,T'}$ the domain bounded by $\Gamma_2,\ \Gamma_3$
and $\Gamma_{4,T'}$  where $\Gamma_{4,T'}$ is the plane 
$\{t=T'\},\ T'\in [\frac{T}{2},T]$.
Let $(u,v)_{\Delta_{1,T'}}$ be the $L_2$ inner product over
$\Delta_{1,T'}$.  Then integrating by parts in
$0=(L_1u,u_t)_{\Delta_{1,T'}}+(u_t,L_1u)_{\Delta_{1,T'}}$
we get,  analogously to (\ref{eq:3.1}), (\ref{eq:3.2}),
with $\Delta_1$ replaced  by $\Delta_{1,T'}$:
\begin{equation}                              \label{eq:3.3}
\|u_t\|_{0,\Gamma_{4,T'}}^2 + \|u\|_{1,\Gamma_{4,T'}}^2\leq
C(\|u\|_{1,\Gamma_{2,T'}}^2 +\|u\|_{1,\Gamma_{3,T'}}^2
+\int_{\frac{T}{2}}^{T'}(\|u_t\|_{0,\Gamma_{4,t}}^2+\|u\|_{1,\Gamma_{4,t}}^2)dt.
\end{equation}
Here 
$\Gamma_{2,T'},\ \Gamma_{3,T'}$ are parts of $\Gamma_2$ and $\Gamma_3$
where $t\leq T'$.
Since $T$ is small we get from (\ref{eq:3.3}) that
\begin{equation}                              \label{eq:3.4}
 \max_{T'\in [\frac{T}{2},T]}(\|u_t\|_{0,\Gamma_{4,T'}}^2 + 
\|u\|_{1,\Gamma_{4,T'}}^2)\leq
C(\|u\|_{1,\Gamma_2}^2 + \|u\|_{1,\Gamma_3}^2).
\end{equation}

For any $w_0,w_1\in C_0^\infty(\Gamma_4)$ there exists a 
smooth solution
of the  Cauchy
 problem $L_1u=0,\ u(y,T)=w_0,\ \frac{\partial u(y,T)}{\partial t}=w_1$
in the domain $t<T$ (see, for example, [H]).  In particular,
$\{u|_{\Gamma_2},u|_{\Gamma_3}\}\in \mathcal{H}_1$. 
  
We shall show that the image of the map $\{w_0,w_1\}\rightarrow 
\{u|_{\Gamma_2},u|_{\Gamma_3}\}$
where $\{w_0,w_1\}\in \mathcal{H}$ are smooth, $L_1u=0$, is dense in $\mathcal{H}_1$.
  Suppose there exists $\{\varphi,\psi\}\in \mathcal{H}_1$ 
such that
$(v,\varphi)_1+(v_1,\psi)_1=0,\ \ \forall \{v,v_1\}\in \mathcal{H}_1,
\ v=u|_{\Gamma_2}, \ v_1=u|_{\Gamma_3}$, $L_1u=0,\ u$ is smooth.

As in the
Lemma \ref{lma:2.3}
(see the remark after the end of the proof of Lemma \ref{lma:2.3})
 one can find $\{\varphi_0,\psi_0\}\in 
\mathcal{H}_1$ such that 
$Q(\varphi_0,v)+Q_1(\psi_0,v_1)=(\varphi,v)_1+(\psi,v)_1$
for any $\{v,v_1\}\in \mathcal{H}_1$.
Therefore $Q(\varphi_0,v)+Q_1(\psi_0,v_1)=0,
\ \forall \{v,v_1\},\ v=u|_{\Gamma_2},\ v_1=u|_{\Gamma_3},
\ u$  is smooth,  $L_1u=0$.
Here $(\varphi,v)_1$ is the inner product in 
$H^1(\Gamma_j),\ j=2,3$.

Extend $\varphi_0(s,y')$ by 0  for
$s>T$ and  
extend 
$\psi_0(\tau,y')$ by 0 for $\tau<-T$.  Let 
\[
b(y',y_n,t)=
\varphi_0(s,y')
+
 \psi_0(\tau,y')
-
\varphi_0(0,y').
\]
  Note that $\varphi_0(0,y')=\psi_0(0,y'),\ b|_{\Gamma_2}=\varphi_0(s,y'),\ 
b|_{\Gamma_3}=\psi_0(\tau,y')$.
Let $\tilde{\Delta}_1$ be the region $\tau<0,s>0$.
Note that $\tilde{\Delta}_1\supset \Delta_1$.
Consider
$L_1^*b=(-4\frac{\partial^2}{\partial s\partial\tau} + (L_1')^*)b$.
We have $\frac{\partial^2 b}{\partial s \partial \tau}=0$ 
in $\tilde{\Delta_1}$ and $(L_1')^*b\in H_{0,-1}$ in $\tilde{\Delta_1}$ 
where $H_{p,r'}$ is the Sobolev space of order $p$ in all variables and
of order $r'$ in $y'$.
Let $f=-L_1^* b$ in $\tilde{\Delta}_1$ and $f=0$ otherwise.
Then $f\in H_{0,-1}(\R^{n+1})$ and $f=0$ for $t<\frac{T}{2}$.
It follows from [E2],  for example,  that there exists 
  $u_0\in H_{1,-1}(\R^{n+1})$
with the weight $e^{-\sigma t}$ such that $L_1^*u_0=f$ in $\R^{n+1},\ 
u_0=0$ for $t<\frac{T}{2}$.  Denote $u^{(0)}=b(y',y_n,t)+u_0$.
Then $L_1^*u^{(0)}=0$ in $\tilde{\Delta}_1,\ u^{(0)}|_{\Gamma_2}=\varphi_0$,
$u^{(0)}|_{\Gamma_3}=\psi_0$.
We  used that $u_0=0$ on $\Gamma_3$  and $\Gamma_2$
by the finite domain of
dependence property since $u_0=0$ for $t<\frac{T}{2}$.

Denote $u_{10}=u^{(0)}|_{t=T},\ u_{20}
=\frac{\partial u^{(0)}}{\partial t}\left|_{t=T}\right. .$ 
Applying Green formula (\ref{eq:3.1}) to $u^{(0)}$ and $v$,  where $v$
is smooth,  $L_1v=0,\ v|_{\Gamma_4}\in C_0^\infty(\Gamma_4)$  we get
\begin{equation}                             \label{eq:3.5}
E(u^{(0)},v)=Q(u^{(0)},v)+Q_1(u^{(0)},v)=0.
\end{equation}
 Since
$u^{(0)}|_{t=T}$ and 
$\frac{\partial u^{(0)}}{\partial y}\left|_{t=T}\right.$
are distributions, the pairing in (\ref{eq:3.5}) is understood as
an extension of the $L_2$-inner product.
Since $v|_{t=T}$ and $\frac{\partial v}{\partial t}\left|_{t=T}\right.$
are arbitrary and smooth we get that
$u^{(0)}|_{t=T}=\frac{\partial u^{(0)}}{\partial t}\left|_{t=T}\right.=0.$
Then by the uniqueness of the Cauchy problem (see [H]) we get that
$u^{(0)}=0$ in $\Delta_1$.  In particular,  $u^{(0)}|_{\Gamma_2} 
=\varphi_0=0,\ u^{(0)}|_{\Gamma_3}=\psi_0=0$.
It follows from (\ref{eq:3.4}) that the image of the map
$\{w_0,w_1\}\in \mathcal{H}
\rightarrow \{u|_{\Gamma_2},u|_{\Gamma_3}\}\in \mathcal{H}_1$ is onto.

Therefore for any $\{v^{(0)},w^{(0)}\}\in \mathcal{H}_1$
there exists $v^{(1)}\in 
\stackrel{\circ}{H^1}(\Gamma_4),\ v^{(2)}\in L_2(\Gamma_4)$
and $v\in H_1(\Delta_1)$ such that 
\begin{equation}                               \label{eq:3.6}
L_1v=0\ \mbox{in}\ \ \Delta_1,\ \ v|_{\Gamma_2}=v^{(0)},\ \ 
\ \ v|_{\Gamma_3}=w^{(0)},\ \ 
v|_{t=T}=v^{(1)},\ \ \frac{\partial v}{\partial t}\left|_{t=T}\right.=
v^{(2)}.
\end{equation}
Lemma \ref{lma:3.1} is proven.
\qed

Denote by $\Delta_2$ the domain bounded by $t=T,\ y_n=0$  and $t-y_n=0$.
Note that $\Delta_1\subset \Delta_2$.  There exists (see, for example, [H])
a unique solution $u\in H^1(\Delta_2)$ of
$L_1u=0$ in $\Delta_2$ with the initial conditions 
 $u|_{t=T}=v^{(1)},\ \frac{\partial u}{\partial t}\left|_{t=T}\right.
=v^{(2)}$
and the boundary condition $u|_{y_n=0}=0$.  Here $v^{(1)}$ and
$v^{(2)}$ are the same as in (\ref{eq:3.6}).  By the uniqueness of
the Cauchy problem $u=v$ in $\Delta_1\subset\Delta_2$, where $v$ is 
the solution obtained in  (\ref{eq:3.6}).

Therefore we proved the following lemma:
\begin{lemma}                               \label{lma:3.2}
For any $v_0\in 
\stackrel{\circ}{H^1}(Y_{20})\subset 
\stackrel{\circ}{H^1}(\Gamma_2)$ there exists 
$u\in H^1(X_{20})$
such that $L_1u=0$ in $X_{20},\ u|_{Y_{20}}=v_0,\ u|_{y_n=0}=0$.
\end{lemma}
Analogous result holds when $X_{20}$ is replaced by $X_{2s_0}$
and $L_1$ is replaced by $L_1^*$.

Finally we can finish the proof of Lemma \ref{lma:2.2}.
Suppose that the set $H$ of $\{u^g\},\ g\in 
C_0^\infty(\Delta_{2s_0})$
is not dense in $\stackrel{\circ}{H^1}(R_{s_0})$.
Then there exists $w_0\in
 \stackrel{\circ}{H^1}(Y_{2s_0})$ such that 
$(u^g,w_0)_1=0,\ \forall u^g,\ g\in C_0^\infty(\Delta_{2s_0})$
and $(w_0,v)=1$ for some $v\in \stackrel{\circ}{H^1}(R_{s_0})$.
By Lemma \ref{lma:2.3}
one can find $w\in \stackrel{\circ}{H^1}(Y_{2s_0})$ 
such that $Q(u^g,w)=(u^g,w_0)_1$
for all $u^g$.  By Lemma \ref{lma:3.2}
there exists $u\in H^1(X_{2s_0})$ such that
$L_1^*u=0$ in $X_{2s_0},\ u|_{y_n=0}=0,\ u|_{Y_{2s_0}}=w.$

Applying the Green's formula (\ref{eq:2.24}) we get
\[
0=Q(u^g,w)=-\int_{\Delta_{2s_0}}\frac{\partial u^g}{\partial y_n}
\overline{u_t}dy'dt -
\int_{\Delta_{2s_0}}\frac{\partial g}{\partial t}
\frac{\overline{\partial u}}{\partial y_n}dy'dt
\]
since $u^g|_{y_n=0}=g$.  
Since $u|_{y_n=0}=0$ we get  
that 
\[
\int_{\Delta_{2s_0}}\frac{\partial g}{\partial t}
\frac{\overline{\partial u}}{\partial y_n}
dy'dt=0, \ \forall g\in C_0^\infty(\Delta_{2s_0}).
\]

Since $g$ is arbitrary we get that $\frac{\partial}{\partial t}
\frac{\partial u}{\partial y_n}=0$ in $\Delta_{2s_0}=\Gamma^{(1)}
\times(s_0,T]$.
Therefore $\frac{\partial u}{\partial t}$ satisfies
$L_1^*\frac{\partial u}{\partial t}=0$ in
$X_{2s_0},\ \frac{\partial u}{\partial t}= 
\frac{\partial }{\partial y_n}\frac{\partial u}{\partial t}=0$ on 
$\Delta_{2s_0}$.  By the unique continuation theorem  (see [T])
we get that $\frac{\partial u}{\partial t}=0$ in 
the double cone of influence of $\Gamma^{(1)}\times[s_0,T]$.
Note that the intersection of this double cone with the plane $\tau=0$
is
$R_{s_0}$.
Therefore $w=u\left|_{Y_{2s_0}}\right.$ satisfies the elliptic
equation $(4\frac{\partial^2}{\partial s^2} +(L_1')^*)w=0$ on $R_{s_0}$
since
$\frac{\partial^2 u}{\partial t^2}=0,\ 
\frac{\partial^2 u}{\partial s^2}=\frac{1}{4}
(\frac{\partial}{\partial t} -\frac{\partial}{\partial y_n})^2u
=\frac{1}{4}\frac{\partial^2 u}{\partial y_n^2}$.

Consider $Q(w,v)$.  Note that $\mbox{supp}\ v\in  \stackrel{\circ}{H^1}(R_{s_0})$.
Integrating by parts and using that $w$ satisfying the equation 
$4w_{ss}+(L_1')^*w=0$ in $R_{s_0}$ we get that $Q(w,v)=0$ and this 
contradicts the condition that $Q(w,v)=1$.
Therefore $H$ is dense in $\stackrel{\circ}{H^1}(R_{s_0})$.
\qed

Let $L^{(p)}u^{(p)} = 0,\ p=1,2,$ be two hyperbolic equations of 
the form (\ref{eq:1.1})
in domains $\Omega^{(p)},\ p=1,2,$ respectively, satisfying the initial-boundary 
conditions (\ref{eq:1.2}). We assume that 
$\Gamma_0\subset\partial\Omega^{(1)}\cap\partial\Omega^{(2)}$, 
$\mbox{\mbox{supp}\ }f\subset \Gamma_0\times(0,T_0]$ and
$\Lambda^{(1)}=\Lambda^{(2)}$ on $\Gamma_0\times (0,T_0)$  where
$\Lambda^{(p)}$ is the D-to-N operator corresponding to $L^{(p)},\ p=1,2.$

Let $B\subset\Omega^{(1)}\cap\Omega^{(2)}$ be homeomorphic to a ball and
such that the domains $\Omega_i=\Omega^{(i)}\setminus\overline{B},i=1,2,$ 
are smooh.
Assume that $\partial B\cap\partial\Omega^{(i)}\subset\Gamma_0,i=1,2,$
and connected.  Let  $\gamma^{(i)}=\partial\Omega^{(i)}\setminus\overline{\Gamma_0},
i=1,2.$  Note that $\partial\Omega_1\setminus\gamma^{(1)}=\partial\Omega_2
\setminus\gamma^{(2)}\stackrel{def}{=}\Gamma_1$.

We conclude this section with the following lemma that will be important
in the global step (c.f. [KKL1],  Lemma 9):

\begin{lemma}                         \label{lma:3.3}
Let $L^{(1)}=L^{(2)}$ in $\overline{B}$ and let
$\delta=\max_{x\in \overline{B}}d(x,\Gamma_0)$,
where $d(x,\Gamma_0)$ is the distance in $\overline{B}$ from
$x \in \overline{B}$ to $\Gamma_0$.
Let $\Lambda_i$  be the D-to-N operators corresponding to $L^{(i)}$
in smaller domains   $\Omega_i,i=1,2$.  If $\Lambda^{(1)}     
=\Lambda^{(2)}$  on $\Gamma_0\times(0,T_0)$  then $\Lambda_1=\Lambda_2$
on $\Gamma_1\times(\delta,T_0-\delta)$.
\end{lemma}

{\bf Proof:}
Let $\Omega$ be either $\Omega^{(1)}$ or $\Omega^{(2)}$.
Denote by $\Omega_T=\Omega\times(-\infty,T),\ 
\partial\Omega_T=\partial\Omega\times (-\infty,T).$  
Let $H_s^+(\Omega_T)$ be the Sobolev space with norm $\|u\|_s$
where $u=0$ for $t<0$.
The proof of Lemma \ref{lma:3.3}  will be based on 
a  variant of the Runge theorem.  
\begin{lemma}                             \label{lma:3.4}
Denote by $\tilde{\Omega}=
\Omega\setminus \overline{B},\ \gamma=\partial\Omega\setminus\overline{\Gamma_0}.$
Let $u\in H_{1+s}^+(\tilde{\Omega}_{T})$ be the solution of
\begin{equation}                                \label{eq:3.7}
\begin{array}{l}
Pu=0,\ \ (x,t)\in \tilde{\Omega}_{T},
\\
u|_{\partial\tilde{\Omega}_{T}}=f,\ \ \ u=0\ \ \mbox{for}\ \ t<\delta,
\end{array}
\end{equation}
where $P$ has the form (\ref{eq:1.1}), 
$f\in  H_{s+1}^+(\partial\tilde{\Omega}_{T}),\ f=0$ on 
$\gamma_T=\gamma\times(-\infty,T),\ f=0$ for $t<\delta,\  s\leq 0$.
Then there exists a sequence of smooth functions $u_n$ in $\Omega_T$  such that
\begin{equation}                                 \label{eq:3.8}
Pu_n=0,\ \ (x,t)\in \Omega_{T}, 
\ \
u_n|_{\partial\Omega_{T}}=f_n, 
\end{equation}
where $u_n=0$ for $t<0$, $f_n=0$ on $\gamma_T$ and
$u_n\rightarrow u$ when $n\rightarrow \infty$
in $H_{s}^+(\tilde{\Omega}_{T})$.  
\end{lemma}
{\bf Proof:} 
Let $\Omega_\e$ be a smooth 
domain in $\R^n$ such that $\Omega_\e\supset\Omega,
\ \partial\Omega_\e=\gamma\cup\gamma_\e,$ 
and $\Gamma_0$ is inside of $\Omega_\e$.

The following existence and uniquenesss lemma is well known (see, for
example, [E2]).
\begin{lemma}                             \label{lma:3.5}
For any $h\in H_{s}^+(\Omega_T)$ and $f\in H_{s+1}^+(\partial\Omega_T),
\ s\geq 0$, there exists a unique $u\in H_{1+s}^+(\Omega_T)$ 
such that 
\[
Pu=h\ \ \mbox{in}\ \ \Omega_T,\ \ u|_{\partial\Omega_T}=f.
\]
Moreover, $\frac{\partial u}{\partial \nu}\left|_{\partial\Omega_T}
\in H_s^+(\partial\Omega_T)\right.$ and
\[
\|u\|_{1+s}+\|\frac{\partial u}{\partial \nu}\|_s\leq
C\|f\|_{s+1}+C\|h\|_{s}. 
\]
When $h=0$ the same result is true for any $s\in\R$.
\end{lemma}
The last statement follows from the fact that when $Pu=0$ and
$s<0$ the norm $\|u\|_{1+s}$ is equivalent near 
$\partial\Omega_T$
 to the norm $\|u\|_{1,s}$ (c.f. [E2]).
Here $\|u\|_{p,s}$ is the Sobolev norm of order $p$ in all variables
and of order $s$ in $(y',t)$.

Denote by $E(x,t,y,\tau)$ the forward Green function for the domain 
$\Omega_\e$.  More precisely  $E(x,t,y,\tau)$ satisfies
\[
\begin{array}{l}
PE=\delta(x-y)\delta(t-\tau),\ (x,t)\in\Omega_{\e}\times(-\infty,\infty),\ 
(y,\tau)\in \Omega_\e\times(-\infty,+\infty),
\\
E=0\ \ \mbox{for}\ \ t<\tau,\\  E=0\ \ \mbox{when}\ \ (x,t)\in
\partial\Omega_\e\times(-\infty,+\infty),
\ \ (y,\tau)\in \Omega_\e\times(-\infty,+\infty).
\end{array}
\]
The existence of $E(x,t,y,\tau)$ follows from the Hadamard construction
(see [H]) and the Lemma \ref{lma:3.5}.

Denote by $E$ the operator with the kernel $E(x,t,y,\tau)$.
Let  $D_\e =(\Omega_\e\setminus\overline{\Omega})\times(0,T+\delta)$.
Denote by $K$ the closure  of 
$E\varphi,\ \forall\varphi\in C_0^\infty(D_\e)$, in the norm
of $H_{s}^+(\tilde{\Omega}_{T}),\  s\leq 0$.
Let $K^\perp\subset \stackrel{\circ}{H_{-s}^-}(\tilde{\Omega}\times(0,T))$ be 
the "orthogonal" complement to $K$,  i.e.
$(f,g)=0$ for all $f\in K \ \ \mbox{iff}\ \ g\in K^\perp$.  Here
$(f,g)$ is the extension of the $L_2$ scalar product in
$\tilde{\Omega}\times (0,T)$ and 
$\stackrel{\circ}{H_{-s}^-}(\tilde{\Omega}\times(0,T))$ 
consists of functions 
that belong to $H_{-s}(\R^n\times(0,+\infty))$ after
being extended by zero for $t>T$ and for $x\not\in \overline{\tilde{\Omega}}$.
 Note that $(f,g)=(lf,g_+)$ where 
$lf\in H_{s}(\R^n\times(-\infty,+\infty))$  and
$g_+\in H_{-s}(\R^n\times (-\infty,+\infty))$ are arbitrary
extensions of $f$ and $g$  such that $g_+=0$  for 
$x\not\in\overline{\tilde{\Omega}},
g_+=0$  for $t> T, lf=0$  for $t<0$.
Let $g\in K^\perp$ be arbitrary and let $\overline{w(y,\tau)}=
\overline{E^*g_+}$ or, formally:
\[
\overline{w(y,\tau)}=
\int_{-\infty}^\infty\int_{\R^n} E(x,t,y,\tau)
\overline{g_+(x,t)}dxdt.
\]
We have $(E\varphi,g_+)=(\varphi,E^*g_+)=0$ for any
$\varphi\in C_0^\infty(D_\e)$.
Therefore $w(y,\tau)=0$ for 
$(y,\tau)\in D_\e,\ 
w(y,\tau)$ is the solution of the initial-boundary value problem
\begin{equation}                                     \label{eq:3.9}
\begin{array}{l}
P^*w=g_+(y,\tau),\ \ (y,\tau)\in (\Omega_\e\times (0,+\infty)),
\\
w(y,\tau)=0\ \ \mbox{for}\ \ \tau>T,\ \ y\in \Omega_\e,
\ \ w|_{\partial\Omega_\e\times(0,+\infty)}=0,
\end{array}
\end{equation}
and $\overline{E(x,t,y,\tau)}$ is the Schwartz kernel of 
the solution operator to (\ref{eq:3.9}). Here $P^*$ is the adjoint to
$P$. 

It follows from Lemma \ref{lma:3.5}   that 
$w\in H_{1-s}(\Omega_\e\times(0,+\infty)),
\ w=0$ for $t>T,
\ w=0$ on 
$\partial\Omega_\e\times(0,+\infty)$ and
$\frac{\partial w}{\partial \nu}
\left|_{\partial\Omega_\e\times(0,+\infty)}\right.
\in H_{-s}(\partial\Omega_\e\times(0,+\infty))$.   
Since $w(y,\tau)=0$ in $D_\e$
 and $g_+=0$ in $(\Omega_\e\setminus\overline{\tilde{\Omega}})
\times(0,+\infty)$ we get by the uniqueness of the Cauchy problem 
(see [T]) that $w=0$ in $(\Omega_\e\setminus\overline{\tilde{\Omega}})
\times(\delta,+\infty)$.   Therefore $w$ and $\frac{\partial w}{\partial \nu}$
are zero on $(\partial\tilde{\Omega}\setminus\gamma)\times(\delta,+\infty)$.
Take any
$u\in H_{1+s}^+(\tilde{\Omega}_{T})$ that satisfies (\ref{eq:3.7}).
We have $(u,g_+)=(u,P^*w)$.  Using the Green's formula 
in $\tilde{\Omega}\times (-\infty,\infty)$
we get
$(u,g_+)=(Pu,w)=0$, i.e. $u\in K$.
\qed

Now we shall prove Lemma \ref{lma:3.3}.

Let  $v_i,\ i=1,2,$ be smooth solutions of $L^{(i)}v_i=0$ in 
$\Omega_{i}\times(-\infty,T_0-\delta),\ v_i=0$ for  $t<\delta$
and such that $v_i=f$ on $\Gamma_1\times(-\infty,T_0-\delta)
,\ \mbox{\mbox{supp}\ }f\subset \Gamma_1\times(-\infty,T_0-\delta)$.
Extend $f$ smoothly from $\partial\Omega_i\times(-\infty,T_0-\delta)$
to
$\partial\Omega_{i}\times(-\infty,T_0),\ \mbox{\mbox{supp}\ }f\subset\Gamma_1\times(-\infty,T_0)$.  
Then extend $v_i$ in $\Omega_{i,T_0}$ as solutions of
$L^{(i)}v_i=0,\ v_i=f$ on $\partial\Omega_{i,T_0}$.
  We shall show that 
$\Lambda_1 f=\Lambda_2 f$ on $\Gamma_1\times (\delta,T_0-\delta)$.
By Lemma \ref{lma:3.4} there exists a sequence 
of smooth functions $w_{n1}$ in $\Omega^{(1)}_{T_0}$
 satisfying $L^{(1)}w_{n1}=0$ such that 
$\|v_1-w_{n1}\|_{s}\rightarrow 0$ when $n\rightarrow \infty$.  Here
$\|\  \|_{s}$ is the norm in 
 $H_{s}(\Omega_{1,T_0})$.
Denote $f_n=w_{n1}|_{\partial\Omega_{T_0}^{(1)}}\in 
H_{s+\frac{1}{2}}^+(\partial\Omega_{T_0}^{(1)})$.
  Note
that $\mbox{\mbox{supp}\ }f_n\subset \Gamma_{0,T_0}$.
By Lemma \ref{lma:3.5} there exists $w_{n2}\in 
H_{s+\frac{1}{2}}^+(\Omega_{T_0}^{(2)})$  
such that $L^{(2)}w_{n2}=0$ in
$\Omega_{T_0}^{(2)},\ w_{n2}=f_n$ on $\partial\Omega_{T_0}^{(2)},
\ f_n=0$ on $\gamma_{T_0}^{(2)}$.  
Since $\Lambda^{(1)}=\Lambda^{(2)}$ on $\Gamma_{0,T_0}$ we have
$\Lambda^{(1)}f_n=\Lambda^{(2)}f_n$ on $\Gamma_{0,T_0}$.  Since $L^{(1)}=L^{(2)}$
in $B$  
the uniqueness of the Cauchy problem (see [T]) implies that
$w_{n1}=w_{n2} $ 
and $\frac{\partial w_{n1}}{\partial \nu}=
\frac{\partial w_{n2}}{\partial \nu}$ on 
$\Gamma_1\times(-\infty,T_0-\delta)$.

It follows from $L^{(1)}w_{n1}=0$ and 
$w_{n1}\in H_{s}^+(\Omega_{1,T_0})$
that $ w_{n1}\in
 H_{2,s-2}^+$  near $\partial\Omega_{1,T_0}$ (c.f. [E2]).  Therefore
 $w_{n1}|_{\partial\Omega_{1,T_0}}=f_{n1}\in 
H_{s-\frac{1}{2}}(\partial\Omega_{1,T_0}), 
\frac{\partial w_{n1}}{\partial \nu}\left|_{\partial\Omega_{1,T_0}}\right.
=g_{n1}\in H_{s-\frac{3}{2}}^+(\partial\Omega_{1,T})$ 
and 
$\|f_{n1}-f\|_{s-\frac{1}{2}}+\|g_{n1}-
\frac{\partial v_1}{\partial \nu}\|_{s-\frac{3}{2}}
\leq C\|w_{n1}-v_1\|_{s}\rightarrow 0$
when $n\rightarrow \infty$.  Denote 
\[
\frac{\partial v_i}{\partial \nu}\left.\right|_{\partial\Omega_{i}
\times(\delta,T-\delta)}=g_i,\ i=1,2,
\ f_{n2}=w_{n2}\left.\right|_{\partial\Omega_2\times(0,T-\delta)},
\ g_{n2}= \frac{\partial w_{n2}}{\partial \nu}\left.\right|_{\partial\Omega_{2}\times
(0,T-\delta)}.
\] 
Since $f_{n1}=f_{n2}$ on $\Gamma_1\times(0,T-\delta)$ we have that 
$\|f_{n2}-f\|_{s-\frac{1}{2}}'\rightarrow 0$ when $n\rightarrow\infty$
where $\|\ \|_s'$ means the norm in 
$H_s^+(\partial\Omega_{2,T-\delta})$.
Therefore applying Lemma \ref{lma:3.5} to $v_2-w_{n2}$ we get that
$\|g_2-g_{n2}\|_{s-\frac{3}{2}}'\leq C\|f_{n2}-f\|_{s-\frac{1}{2}}'\rightarrow 0$
when $n\rightarrow \infty$.  Since $g_{n1}=g_{n2}$ on 
$\Gamma_1 \times(0,T-\delta)$ we get that
$g_1=g_2$ on $\Gamma_1\times(\delta,T-\delta)$.
\qed
\\
\
\\

{\bf Acknowledgements.}
Author expresses deep gratitude to Jim Ralston 
for his generous help.  My special thanks  to Slava Kurylev 
who saved me from embarrassing mistakes and gave many useful suggestions
that greatly improved this paper.

\end{document}